\newfont{\Bbb}{msbm10 scaled\magstephalf}
\documentclass[12pt]{amsart}
\usepackage{amsmath, amsthm, amscd, amsfonts}
\newtheorem{thm}{Theorem}
 
 \newtheorem{lem}[thm]{Lemma}
 
 \def\be{\begin{equation}}
 \def\ee{\end{equation}}
 \def\bea{\begin{equation}\begin{array}{rcl}}
 \def\eea{\end{array}\end{equation}}

\begin{document}
\title[EXTENDED CES$\acute{A}$RO OPERATORS]{{EXTENDED CES$\acute{A}$RO OPERATORS ON ZYGMUND SPACES IN THE UNIT BALL} }
\author[Z.S. Fang and Z.H.Zhou]{Zhong-Shan Fang \and Ze-Hua Zhou$^*$ }
\address{\newline Department of Mathematics\newline
Tianjin Polytechnic University
\newline Tianjin 300160\newline P.R. China.}

\email{fangzhongshan@yahoo.com.cn}

\address{\newline Department of Mathematics\newline
Tianjin University
\newline Tianjin 300072\newline P.R. China.}
\email{zehuazhou2003@yahoo.com.cn}



\date{}
\thanks{\noindent $^*$ Ze-Hua Zhou, Corresponding author. Supported in part by the National Natural Science Foundation of
China (Grand Nos.10671141, 10371091).}

\begin{abstract}
Let $g$ be a holomorphic function of the unit ball $B$ in the
$n$-dimensional space, and denote by $T_g$ and $I_g$ the induced
extended Ces$\acute{a}$ro operator and another integral operator.
The boundedness and compactness of $T_g$ and $I_g$ acting on the
Zygmund spaces in the unit ball are discussed and necessary and
sufficient conditions are given in this paper.
\end{abstract}

\maketitle
\section{Introduction}

Let $f(z)$ be a holomorphic function  on the unit disc $D$ with
$Taylor$ expansion $f(z)=\sum\limits^{\infty}_{j=0}a_jz^j$, the
classical Ces$\acute{a}$ro operator acting on $f$ is
$${\mathcal C}[f](z)=\sum\limits^{\infty}_{j=0}\left(\frac{1}{j+1}\sum\limits^j_{k=0}a_k\right)z^j.$$

In the past few years, many authors focused on the boundedness and
compactness of extended Ces$\acute{a}$ro operator between several
spaces of holomorphic functions . It is well known that the operator
${\mathcal C}$ is bounded on the usual Hardy spaces $H^p(D)$ for
$0<p< \infty$ and Bergman space, we recommend the interested readers
refer to \cite{Si1, Si3, Mia, DS, Si4}. But the operator ${\mathcal
C}$ is not always bounded, in \cite{SR}, Shi and Ren  gave a
sufficient and necessary condition for the operator ${\mathcal C}$
to be bounded on mixed norm spaces in the unit disc. Recently,
Siskakis and Zhao in \cite{Sis} obtained sufficient and necessary
conditions for Volterra type operator, which is a generalization of
${\mathcal C}$ , to be bounded or compact between $BMOA$ spaces in
the unit disc. It is a natural question to ask what are the
conditions for higher dimensional case.

Let $dv$ be the $Lebesgue$ measure on the unit ball $B$ of $C^n$
normalized so that $v(B)=1$,and $dv_{\beta}=c_\beta(1-|z|^2)^\beta
dv$,where $c_\beta$ is a normalizing constant so that $dv_\beta$ is
a probability measure. The class of all holomorphic functions on $B$
is defined by $H(B)$ . For $f\in H(B)$ we write
$$
 Rf(z)=\sum\limits^n_{j=1}z_j\frac{\partial f}{\partial
z_j}(z).$$

A little calculation shows ${\mathcal
C}[f](z)=\frac{1}{z}\int^z_0f(t)(\log\frac{1}{1-t})'dt$. From this
point of view, if $g\in H(B)$, it is natural to consider the
extended Ces$\acute{a}$ro operator (also called Volterra-type
operator or Riemann-Stieltijes type operator) $T_g$ on $H(B)$
defined by
$$T_g(f)(z)=\int^1_0f(tz)Rg(tz)\frac{dt}{t}.$$
It is easy to show that $T_g$ take $H(B)$ into itself. In general,
there is no easy way to determine when an extended Ces$\acute{a}$ro
operator is bounded or compact.

Motivated by \cite{SR}, Hu and Zhang \cite {Hu1,Hu2,Zha} gave some
sufficient and necessary conditions for the extended ${\mathcal C}$
to be bounded and compact on mixed norm spaces, Bloch space as well
as Dirichlet space in the unit ball.

Another natural integral operator is defined as follows:
$$
I_g(f)(z)=\int^1_0Rf(tz)g(tz)\frac{dt}{t}.$$ The importance of them
comes from the fact that
\begin{equation}T_g(f)+I_g(f)=M_gf-f(0)g(0)\label{1}\end{equation} where the
multiplication operator is defined by
$$M_g(f)(z)=g(z)f(z) , f\in H(B), z\in B.$$

Now we introduce some spaces first. Let $H^{\infty}$ denote the
space of all bounded holomorphic functions on the unit ball,
equipped with the norm $||f||_\infty=\sup\limits_{z\in B}|f(z)|.$

The Bloch space $\mathcal{B}$ is defined as the space of holomorphic
functions such that
$$||f||_\mathcal{B}=\sup\{(1-|z|^2)|Rf(z)|:z\in B\}<\infty.$$
It is easy to check that if $f\in \mathcal{B}$ then
\begin{equation}|f(z)|\leq C\log\frac{2}{1-|z|^2}\|f\|_\mathcal{B}.\label{2}\end{equation}

We define weighted Bloch space $\mathcal{B}_{log}$ as the space of
holomorphic functions $f\in H(B)$ such that
$$||f||_{\mathcal{B}_{log}}=\sup\{(1-|z|^2)|Rf(z)|\log\frac{2}{1-|z|^2}:z\in B\}<\infty.$$

The Zygmund space  $\mathcal{Z}$ \cite{Zhu1} in the unit ball
consists of those functions whose first order partial derivatives
are in the Bloch space.

It is well known that (Theorem 7.11 in \cite{Zhu1}) $f\in
\mathcal{Z}$ if and only if $Rf\in \mathcal{B}$, and $\mathcal{Z}$
is a Banach space with the norm
\begin{equation}||f||=|f(0)|+||Rf||_\mathcal{B}.\label{3}\end{equation}

The purpose of this paper is to discuss the boundedness and
compactness of extended Ces$\acute{a}$ro operator $T_g$ and another
integral operator $I_g$ on the Zygmund space in the unit ball.

\section{Some Lemmas}

In the following, we will use the symbol $C$ to denote a finite
positive number which does not depend on variable $z$ and $f$.

In order to prove the main results, we will give some Lemmas first.

\begin{lem}Assume $f\in \mathcal{Z}$, then we have
$$||f||_{\infty}\leq C||f||$$
\end{lem}
{\bf Proof.} Since $f\in \mathcal{Z}$ implies that $Rf \in
\mathcal{B}$, it follows from (\ref{2}) that
\begin{equation}|Rf(z)|\leq C
\log\frac{2}{1-|z|^2}\|Rf\|_\mathcal{B}\leq C
\log\frac{2}{1-|z|^2}||f||.\label{4}\end{equation} Furthermore by
$\lim\limits_{|z|\to 1}(1-|z|^2)\log\frac{2}{1-|z|^2}=0$ we have
\begin{equation}(1-|z|^2)|Rf(z)|\leq
C(1-|z|^2)\log\frac{2}{1-|z|^2}||f||<\infty\label{5},\end{equation}
so $f\in \mathcal{B}$. It follows from Theorem 2.2 in \cite{Zhu1}
that
$$Rf(z)=\int_{B}\frac{Rf(z)dv_{\beta}(w)}{(1-<z,w>)^{n+1+\beta}}$$
where $\beta$ is a sufficiently large positive constant. Since
$Rf(0)=0$,
\begin{eqnarray*}
f(z)-f(0)&=&\int_0^1\frac{Rf(tz)}{t}dt
=\int_BRf(w)L(z,w)dv_{\beta}(w)\\
\end{eqnarray*}
where the kernel
$$L(z,w)=\int_0^1(\frac{1}{(1-t<z,w>)^{n+1+\beta}}-1)\frac{dt}{t}$$
satisfies $$|L(z,w)|\leq \frac{C}{|1-<z,w>|^{n+\beta}}$$ for all $z$
and $w$ in $B$. Note that $t^{1/2}\log \frac{2}{t} \leq
\frac{2}{e}\cdot(1-\log 2)$ for all $t\in (0,1]$, then
\begin{eqnarray*}
|f(z)-f(0)|&=&C\int_B\frac{(1-|w|^2)|Rf(w)|dv_{\beta-1}(w)}{|1-<z,w>|^{n+\beta}}\\
&\leq&C\int_B\frac{(1-|w|^2)\log \frac{2}{1-|w|^2}||f|| dv_{\beta-1}(w)}{|1-<z,w>|^{n+\beta}} \\
&\leq&C\int_B\frac{(1-|w|^2)^{1-1/2} ||f||
dv_{\beta-1}(w)}{|1-<z,w>|^{n+\beta}}\\
&\leq&C||f||.
\end{eqnarray*}
The last inequality holds since $\int_B\frac{(1-|w|^2)^t
dv(w)}{|1-<z,w>|^{n+1+t+c}}$ is bounded for $c<0$. This completes
the proof of Lemma 1.

By Lemma 1, Montel theorem and the definition of compact operator,
the following lemma follows.
\begin{lem} Assume that $g\in H(B)$. Then $T_g$ (or $I_g):\mathcal{Z}\rightarrow
\mathcal{Z}$ is compact if and only if $T_g$ (or $I_g)$ is bounded
and for any bounded sequence $(f_k)_{k\in N}$ in $\mathcal{Z}$ which
converges to zero uniformly on $\overline{B}$ as $k\rightarrow
\infty$, $||T_gf_k|| \rightarrow 0$(or $||I_gf_k|| \rightarrow 0$)
as $k\rightarrow \infty.$
\end{lem}

\begin{lem} If $(f_k)_{k\in N}$ is a bounded sequence in $\mathcal{Z}$ which
converges to zero uniformly on compact subsets of $B$ as
$k\rightarrow \infty$, then
$\lim\limits_{k\rightarrow\infty}\sup\limits_{z\in
B}|f_k(z)|=0$.\end{lem}

{\bf proof.} Assume $||f_k||\leq M$. For any given $\epsilon
>0$, there exists $0<\eta<1$ such that
$\frac{\sqrt{1-\eta}}{\eta}<\epsilon.$ Note that $t^{1/2}\log
\frac{2}{t} \leq \frac{2}{e}\cdot(1-\log 2)$ for all $t\in (0,1]$,
then when $\eta<|z|<1$, it follows from (\ref{4}) that
\begin{eqnarray*}
|f_k(z)-f_k(\frac{\eta}{|z|}z)|&=&\left|\int_{\frac{\eta}{|z|}}^1
Rf_k(tz)\frac{dt}{t}\right|\leq C \int_{\frac{\eta}{|z|}}^1
\log\frac{2}{1-|tz|^2}||f_k||\frac{dt}{t} \\
&\leq&
C\frac{|z|}{\eta}\int_{\frac{\eta}{|z|}}^1\frac{||f_k||dt}{(1-|tz|^2)^{1/2}}
\leq C\frac{M}{\eta}\int_{\frac{\eta}{|z|}}^1
\frac{|z|dt}{(1-t|z|)^{1/2}}\\
&\leq&2CM\frac{(1-\eta)^{1/2}}{\eta}<C\epsilon.
\end{eqnarray*}
So we get $\sup\limits_{\eta<|z|<1}|f_k(z)|\leq C\epsilon
+\sup\limits_{|w|=\eta}|f_k(w)|.$ Thus, we have
$$\lim\limits_{k\rightarrow\infty}\sup\limits_{z\in
B}|f_k(z)|\leq\lim\limits_{k\rightarrow\infty}(\sup\limits_{|z|\leq
\eta}|f_k(z)|+\sup\limits_{\eta<|z|<1}|f_k(z)|)\leq C\epsilon.$$ Now
we finish the proof of this lemma.

\begin{lem} Let $g\in H(B)$, then $$
R[T_gf](z)=f(z)Rg(z)$$ for any $f\in H(B)$ and $z\in B$.
\end{lem}

{\bf Proof.}  Suppose the holomorphic function $fRg$ has the
$Taylor$ expansion $$(f R g)(z)=\sum_{|\alpha|\geq
1}a_{\alpha}z^{\alpha}.$$ Then we have
\begin{eqnarray*}
&&R(T_gf)(z)=R\int^1_0 f(tz) R(tz) \frac{dt}{t}
=R\int^1_0\sum_{|\alpha|\geq 1}a_{\alpha}(tz)^{\alpha}\frac{dt}{t}\\
&&=R [\sum_{|\alpha|\geq 1}\frac{a_{\alpha}z^{\alpha}}{|\alpha|}]
=\sum_{|\alpha|\geq 1}a_{\alpha}z^{\alpha}=(f R g)(z).
\end{eqnarray*}

\section{Main Theorems}

\textbf{Theorem 1.} Suppose $g\in H(B)$, then the following
conditions are all equivalent:

(a) $T_g$ is bounded on $\mathcal{Z}$;

(b) $T_g$ is compact on $\mathcal{Z}$;

(c) $g\in \mathcal{Z}.$

{\bf Proof.} $b \Longrightarrow a$ is obvious. For $a\Longrightarrow
c$ we just take the test function given by $f(z)\equiv1$.

We are going to prove $c\Longrightarrow b$. Now assume that $g\in
\mathcal{Z}$ and that $(f_k)_{k\in N}$ is a sequence in
$\mathcal{Z}$ such that $\sup_{k\in N}||f_k||\leq M$ and that $f_k
\rightarrow 0$ uniformly on on $\overline{B}$ as $k\rightarrow
\infty$. Now note that $T_gg_k(0)=0$ and for every $\epsilon
>0$, there is a $\delta \in (0,1)$, such that
$$(1-|z|^2)(\ln\frac{2}{1-|z|^2})^2<\epsilon$$ whenever
$\delta<|z|<1$. Let $K=\{z\in B:|z|\leq \delta\}$, it follows from
Lemma 4 and (4) that
\begin{eqnarray*}
||T_gf_k||&=&\sup\limits_{z\in
B}(1-|z|^2)\left|R(R(T_gf_k))\right|\\
&=&\sup\limits_{z\in B}(1-|z|^2)|Rf_k\cdot Rg+f_k\cdot
R(Rg)|\\
&\leq& \sup\limits_{z\in B}(1-|z|^2)(|Rf_k\cdot Rg|+|f_k\cdot
R(Rg)|)\\
&\leq&\sup\limits_{z\in K}(1-|z|^2)|Rf_k\cdot Rg|+\sup\limits_{z\in
B-K}(1-|z|^2)(|Rf_k\cdot Rg|\\
&&\hspace*{2mm}+\sup\limits_{z\in B}(1-|z|^2)|f_k\cdot
R(Rg)|\\
&\leq&C||g||\sup\limits_{z\in
K}(1-|z|^2)|Rf_k(z)|\log\frac{2}{1-|z|^2}\\
&&\hspace*{2mm}+C||f_k||\cdot||g||\sup\limits_{z\in
B-K}(1-|z|^2)(\log\frac{2}{1-|z|^2})^2+||g||\cdot\sup\limits_{z\in
B}|f_k(z)|.
\end{eqnarray*}
With the uniform convergence of $f_k$ to $0$ and the Cauchy
estimate, the conclusion follows by letting $k\rightarrow \infty$.

\textbf{Theorem 2.} Suppose $g\in H(B)$, $I_g :\mathcal{Z}
\rightarrow \mathcal{Z}$. Then $I_g$ is bounded if and only if $g\in
H^{\infty}\cap\mathcal{B}_{log}$.

{\bf Proof.} First we assume that $g\in
H^{\infty}\cap\mathcal{B}_{log}$. Notice that $I_gf(0)=0$ and
$R(I_gf)=fRg$, it follows from (4) that
\begin{eqnarray*}
(1-|z|^2)|RR(I_g f)(z)|&=&(1-|z|^2)|R(Rf(z)\cdot g(z))|\\
&=&(1-|z|^2)|R(Rf)(z)\cdot g(z)+Rf(z)\cdot Rg(z)|\\
&\leq& \|Rf(z)\|_\mathcal{B}\|g\|_{\infty}+|Rf(z)|(1-|z^2|)|Rg(z)|\\
&\leq& C||f||\cdot ||g||_\infty+C||f||(1-|z|^2)|
Rg(z)|\log\frac{2}{1-|z|^2}\\
&\leq& C||f||\cdot
||g||_\infty+C||f||\cdot||g||_{\mathcal{B}_{log}}.
\end{eqnarray*}
The boundedness of $I_g$ follows.

Conversely, assume that $I_g$ is bounded, then there is a positive
constant $C$ such that
\begin{equation}||I_gf||\leq C||f||\label{7}\end{equation} for every $f\in \mathcal{Z}$.
Setting
$$h_a(z)=(\log\frac{2}{1-|a|^2})^{-1}(<z,a>-1)[(1+\log\frac{2}{1-<z,a>})^2+1]$$
for $a\in B$ such that $|a|\geq\sqrt{1-2/e}$, then
$$Rh_a(z)=<z,a>(\log\frac{2}{1-<z,a>})^2(\log\frac{2}{1-|a|^2})^{-1}$$
and
$$RRh_a(z)=\{<z,a>(\log\frac{2}{1-<z,a>})^2+\frac{2<z,a>^2}{1-<z,a>}\log\frac{2}{1-<z,a>}\}(\log\frac{2}{1-|a|^2})^{-1}$$
It is easy to check that $M=\sup\limits_{\sqrt{1-2/e}\leq
|a|<1}||h_a||<\infty$. Therefore, we have that
\begin{eqnarray}
\infty&>&\|I_g\|\|h_a\|\geq||I_gh_a||\nonumber\\
&\geq&\sup\limits_{z\in B}(1-|z|^2)|RRh_a(z)\cdot g(z)+Rh_a(z)\cdot
Rg(z)|\nonumber\\
&\geq&(1-|a|^2)|\frac{2|a|^4}{1-|a|^2}g(a)+|a|^2\log\frac{2}{1-|a|^2}g(a)+|a|^2Rg(a)\log\frac{2}{1-|a|^2}|\nonumber\\
&\geq&-\{2|a|^4+|a|^2 \frac{2}{e}(1-\log2)\}|g(a)|+|a|^2(1-|a|^2)|Rg(a)|\log\frac{2}{1-|a|^2}\nonumber\\
&\geq&-(2+\frac{2}{e}(1-\log2))|a|^2+|a|^2(1-|a|^2)|Rg(a)|\log\frac{2}{1-|a|^2}.\label{8}
\end{eqnarray}
Next let
$$f_a(z)=h_a(z)-\int_0^1<z,a>\log\frac{2}{1-t<z,a>}dt$$
then
$$Rf_a(z)=<z,a>\{(\log\frac{2}{1-<z,a>})^2(\log\frac{2}{1-|a|^2})^{-1}-\log\frac{2}{1-<z,a>}\}$$
$$RRf_a(z)=RRh_a(z)-<z,a>\log\frac{2}{1-<z,a>}-\frac{<z,a>^2}{1-<z,a>}$$
and consequently $N=\sup\limits_{\sqrt{1-2/e}\leq
|a|<1}||f_a||<\infty$. Note that $Rf_a(a)=0$ and
$RRf_a(a)=\frac{|a|^4}{1-|a|^2}$, we have
\begin{eqnarray}
\infty&>&||I_g||\cdot||f_a||\geq||I_gf_a||\nonumber\\
&\geq&\sup\limits_{z\in B}(1-|z|^2)|RRf_a(z)\cdot g(z)+Rf_a(z)\cdot
Rg(z)|\nonumber\\
&\geq&(1-|a|^2)|RRf_a(a)g(a)+Rf_a(a)Rg(a)|=|a|^4|g(a)|.\label{9}
\end{eqnarray}
From the maximum modulus theorem, we get $g\in H^\infty$. So it
follows from (\ref{8}) and (\ref{9}) that
\begin{equation}\sup\limits_{\sqrt{1-2/e}\leq
|a|<1}(1-|a|^2)|Rg(a)|\log\frac{2}{1-|a|^2}<\infty.\label{10}\end{equation}
On the other hand, we have
\begin{eqnarray}
& &\sup\limits_{|a|\leq
\sqrt{1-2/e}}(1-|a|^2)|Rg(a)|\log\frac{2}{1-|a|^2} \nonumber\\
&&\leq
\frac{2}{e}\cdot(1-\log 2)\max\limits_{|a|=\sqrt{1-2/e}}|Rg(a)|\nonumber\\
&&\leq\sup\limits_{
\sqrt{1-2/e}\leq|a|<1}(1-|a|^2)|Rg(a)|\log\frac{2}{1-|a|^2}<+\infty.\label{11}
\end{eqnarray}
Combining (\ref{10}) and (\ref{11}), we finish the proof of Theorem
2.

\textbf{Corollary } The multiplication operator
$M_g:\mathcal{Z}\rightarrow \mathcal{Z}$ is bounded if and only if
$g\in \mathcal{Z}$.

{\bf Proof.} If $M_g$ is bounded on $\mathcal{Z}$, then setting the
test function $f\equiv1$, we have $M_g f=g \in \mathcal{Z}$.

Conversely, if $g\in \mathcal{Z}$, from Lemma 1 and (\ref{5}), it is
easy to see that $g\in H^{\infty}\cap\mathcal{B}_{log}$, so by
Theorems 1 and 2, both $T_g$ and $I_g$ are bounded, it follows from
(\ref{1}) that $M_g$ is also bounded.

\textbf{Theorem 3.} Suppose $g\in H(B)$, $I_g :\mathcal{Z}
\rightarrow \mathcal{Z}$. Then $I_g$ is compact if and only if
$g=0.$

{\bf Proof.} The sufficiency is obvious. We just need to prove the
necessity. Suppose that $I_g$ is compact, for any given sequence
$(z_k)_{k\in N}$ in $B$ such that $|z_k|\rightarrow 1$ as
$k\rightarrow \infty$, if we can show $g(z_k)\rightarrow 0$ as
$k\rightarrow \infty$, then by the maximum modulus theorem we have
$g\equiv0$. In fact, setting
\begin{eqnarray*}f_k(z)=h_{z_k}(z)-(\log\frac{2}{1-|z_k|})^{-2}\int_0^1 <z,z_k>(\log
\frac{2}{1-t<z,z_k>})^3 dt.\end{eqnarray*} Using the same way as in
Theorem 2, we can show $\sup_{k\in N}||f_k||\leq C$ and $f_k$
converges to $0$ uniformly on compact subsets of $B$. Since $I_g$ is
compact, we have $||I_gf_k||\rightarrow 0$ as $k\rightarrow \infty$.
Note that $Rf_k(z_k)=0$ and $RRf_k(z_k)=-\frac{|z_k|^4}{1-|z_k|^2},$
it follows that
\begin{eqnarray*}
|z_k|^4|g(z_k)|&\leq& \sup\limits_{z\in B}(1-|z|^2)|RRf_k(z)\cdot
g(z)+Rf_k(z)\cdot Rg(z)|\\
&\leq& \sup\limits_{z\in B}(1-|z|^2)|RR (I_gf_k)(z)| \leq
||I_gf_k||\rightarrow 0
\end{eqnarray*}
as $k\rightarrow \infty$. This ends the proof of Theorem 3.

\end{document}